# Logical Foundations of Local Gauge Symmetry and Symmetry Breaking

Yingrui Yang

(yangyri@rpi.edu)

Department of Cognitive Science

Rensselaer Polytechnic Institute, USA

*Abstract*

The present paper reports two results. It is shown that the formula

$$P(x) = \forall y \forall z [\neg G(x,y) \to \neg M(z)]$$

provides the logic underlying gauge symmetry, where $M$ denotes the predicate of being massive. For the logic of spontaneous symmetry breaking, we have

$$P(x) = \forall y \forall z [D(x,y) \to M(z)]$$

Notice that the above two formulas are not logically equivalent. The results are obtained by integrating four components, namely, gauge symmetry and Higgs mechanism in quantum field theory, and Gödel's incompleteness theorem and Tarski's indefinability theorem in mathematical logic. Gödel numbering is the key for arithmetic modeling applied in this paper.

*1. Introduction*

The standard model of particle physics is a well-developed and well-documented theoretical framework. Gauge field theory is the main language adopted in this standard model [1]. Gauge symmetry is an essential systematic property of gauge theoretic modeling, and spontaneous symmetry breaking by Higgs mechanism is the key methodology to supplement the standard model. The standard model has successfully reached some general systematic conclusions. Nonetheless, the metalogic underlying these conclusions has not been spelled out and still remains unclear.

This paper aims to spell out the meta-logical properties underlying the systematic conclusions reached by the standard model of particle physics. Two results will be reported. The first result is about the underlying metalogic of gauge symmetry, showing that the notion of massiveness is indefinable. This result is obtained by integrating the techniques used by Gödel's incompleteness theorem [3] and Tarski's indefinability theorem [5]. The second result is about spontaneous symmetry breaking. Gödel numbering will be introduced, as it is the key to the first-order characterization of Lagrangian. The necessary backgrounds of Gödel's theorem and Tarski's theorem are given in the appendix.

*2. Gödel numbering and first-order characterization*

Mathematical language always deals with symbols, formulas, and derivations. For a mathematical framework, even though its base domains (such as real or complex fields) are uncountable infinities (i.e., the continuum), the number of symbols used to denote variables, functions, operators, etc., is infinite but countably many. Thus, we can have an effective procedure to mechanically assign a unique odd number to each and every symbol in order, called *Gödel number*. For a given symbol *e*, its Gödel number can be written as $g(e)$, which can be seen as a function or an odd number. A formula is a finite string of symbols, written as

$$L = e_1 e_2 \ldots \ldots e_n$$

The Gödel number of a formula can be calculated by

$$g(L) = q_1^{g(e_1)} q_2^{g(e_2)} \ldots \ldots q_n^{g(e_n)}$$

where $q_i$ is the first *i* prime numbers in its natural order, and $g(e_i)$ is the Gödel number of the *i*th symbol in the formula *L*. A derivation is a finite sequence of formulas, written as

$$Der(L) \equiv\, <u_1, u_2, \ldots \ldots, u_m>$$

The Gödel number of a derivation can be calculated by

$$g[(Der(L)] = q_1^{g(u_1)} q_2^{g(u_2)} \ldots q_i^{g(u_i)} \ldots q_m^{g(u_m)}$$

where $g(u_i)$ is the Gödel number of the *i*th formula in the derivation sequence. The Gödel number of any given formula or derivation is always an even number, which is also a composite number.

The above method is called *Gödel numbering* [3]. The beauty and power of Gödel numbering is that, based on the fundamental theorem of arithmetic (i.e., the unique decomposition theorem of primes), from a given Gödel number we can uniquely recapture the original derivation, the original formula, or the original symbol used in the context.

Note that logic has nothing to do with the content. So that the *first-order characterization* of the Lagrangian only requires three conditions. First, the Lagrangian represented by a Gödel number. Second, any given derivation of a particular kind and its Gödel number can be used interchangeably. Third, It allows to introduce new predicates or function terms.

*3. Gauge Particles and the Indefinability of Massiveness*

This section will proceed as follows. First, the general conclusion reached from physics is summarized as a statement. Second, the Indefinability Theorem of Massiveness is stated and the sketch of its proof is shown. Lastly, some comments on this theorem are provided.

*Statement of physics*: The standard model of particle physics is the composite of three gauge symmetry groups: $U(1) \times SU(2) \times SU(3)$. The gauge symmetry implies that all the gauge particles, including gauge fermions and gauge bosons, are massless.

*The Indefinability Theorem of Massiveness*: Under the first-order characterization of the standard model, massiveness is indefinable within its gauge theoretic modeling. In other words, the predicate "being massive", denoted by $M$, has no model; i.e., the model of massive predicate is null.

In the following, we sketch the proof of this theorem. Let $L$ be the Lagrangian of a particle,

$$L = -\frac{1}{4}F_{\mu\nu}F^{\mu\nu} + (D_\mu \varphi)^\dagger D_\mu \varphi - \mu^2(\varphi^\dagger \varphi) + \lambda(\varphi^\dagger \varphi)^2 \tag{1}$$

There is no explicit mass term in this Lagrangian, therefore it is sometimes called a massless Lagrangian. Let the Gödel number of $L$ be $g(L) = i$. Now we assume there was a derivation (a transformation) from $L$ to $L_m$, where $L_m$ is a Lagrangian with mass terms (also known as a massive Lagrangian). Denote such a derivation by $Der(L, L_m)$. This derivation would have a Gödel number, $g[Der(L, L_m)] = j$, and the end formula $L_m$ would have its own Gödel number $g(L_m) = k$.

Now we introduce an arithmetic relation, $\mathcal{G}(i,j)$, which by the definition of expressibility can be expressed by the term $G(\boldsymbol{i}, \boldsymbol{j})$. In addition, we introduce a new predicate of "being massive", denoted by $M$. Then we consider the formula below

$$P(x) = \forall y \forall z [\neg G(x, y) \rightarrow \neg M(z)] \tag{2}$$

Substituting the variables by the corresponding Gödel numbers introduced above, we have

$$S(\boldsymbol{i}, \boldsymbol{j}, \boldsymbol{k}) = \neg G(\boldsymbol{i}, \boldsymbol{j}) \rightarrow \neg M(\boldsymbol{k}) \tag{3}$$

By the statement of physics, there is no such a derivation from $L$ to $L_m$. Hence, $\mathcal{G}(i,j)$ is actually false; since $\mathcal{G}(i,j)$ is expressible by $G(\boldsymbol{i}, \boldsymbol{j})$, we then have $\neg G(\boldsymbol{i}, \boldsymbol{j})$ from the definition of expressibility. Then by propositional logic, we can infer $\neg M(\boldsymbol{k})$. Now assume for contradiction that $M$ had a model:

$$X = \{k \mid k = g(L_m), and\ L_m\ satisfies\ M\} \tag{4}$$

Then, because $L_m$ was supposed to be massive, we would have that $k \in X$; while by $\neg M(\boldsymbol{k})$, we would have $k \notin X$. This is a contradiction. In other words, $X$ must be a null model of $M$, meaning that the massiveness is indefinable within the gauge theoretic language adopted by the standard model of particle physics. QED.

*4. Higgs Mechanism and Definability of Massiveness*

The Higgs mechanism is responsible for causing local gauge symmetry spontaneously to be broken. Because local symmetry involves gauge field, which is massless when the local symmetry is unbroken, spontaneous local gauge symmetry breaking will make a gauge field massive. There are many massless Lagrangians, and we can also make up many massive Lagrangians; which can be connected to which, depend on whether a derivation is available. This should be the meaning of definability for being massive.

Higgs mechanism can be characterized by making four points. First, to force a complex scalar field to be a real field by letting the imaginative part to be zero, and the resulting real scalar field is seen as a ground state, or a vacuum. Second, to introduce another real scalar field, called Higgs field, which causes the fluctuation of the ground state; it is a kind of perturbation. Third, as the result of perturbation, the local symmetry will be broken, and produce a massless field, called Goldstone field. Forth, the Goldstone field will be annihilated, and it somehow returns back to make the gauge field massive. And finally, the mass terms with gauge field usually are eventually constituted by potential parameters and coupling constants. As Ooguri [4] points out, Higgs mechanism does not tell us what is mass; it derives mass as if it exists. The theorem below shows the logic of Higgs mechanism.

$$P(x) = \forall y \forall z [D(x,y) \rightarrow M(z)] \tag{5}$$

Note that this formula is not logically equivalent to (1). We sketch a proof in the following. To restate the notations used in previous section, let the Gödel number of $L$ be $g(L) = i$, $g[Der(L, L_m)] = j$, and $g(L_m) = k$. $D(x,y)$ expresses the binary relation $d(g(L), g(L_m)]$ if $Der(L, L_m)]$ holds. Substituting the variables by the corresponding Gödel numbers introduced above, we have

$$S(\boldsymbol{i}, \boldsymbol{k}, \boldsymbol{j}) = D(\boldsymbol{i}, \boldsymbol{k}) \rightarrow M(\boldsymbol{j}) \tag{6}$$

Note that different from (2), here Godel numbers $\boldsymbol{k}$ and $\boldsymbol{j}$ have switched the places, reflecting our view that the mass is defined by a derivation. Then by propositional logic, we can infer that $M(\boldsymbol{j})$. Now assume that $M$ had a model

$$X = \{j \mid j = g[Der(L, L_m)]\} \tag{7}$$

Many documented results in physics proved that this model is not null. To complete the existential proof, we provide an example below, of which the derivation is partly taken from Aitchison [1]. This example is about Higgs field itself, which is straightforward, and for reader's convenience, it is spelled out step by step without any omission. During the derivation, the key formulas are numbered by (Hi), definitions are numbered by (Dj), and step formulas are not numbered. Notations and definitions are listed below.

$$\varphi = \rho e^{i\theta}, \ \varphi^\dagger = \rho e^{-i\theta} \tag{D1, D2}$$

$$D_\mu \varphi = (\partial_\mu - ieA_\mu)\varphi, \ (D_\mu \varphi)^\dagger = (\partial_\mu + ieA_\mu)\varphi \tag{D3, D4}$$

$$\rho = \tfrac{1}{\sqrt{2}}(v+\chi), \qquad v = \sqrt{\tfrac{\mu^2}{\lambda}} \qquad \text{(D5), (D6)}$$

$$B_\mu \equiv A_\mu - \tfrac{1}{e}\partial_\mu, \qquad M = ve \qquad \text{(D7), (D8)}$$

Consider a massless Langrangian

$$\mathcal{L} = -\tfrac{1}{4}F_{\mu\nu}F^{\mu\nu} + (D_\mu\varphi)^\dagger D_\mu\varphi + \mu^2\varphi^\dagger\varphi - \lambda(\varphi^\dagger\varphi)^2 \qquad \text{(H1)}$$

Since $\varphi^\dagger\varphi = \rho e^{-i\theta}\rho e^{i\theta} = \rho^2 e^0 = \rho^2$

$$D_\mu\varphi = (\partial_\mu - ieA_\mu)\rho e^{i\theta} = \partial_\mu(\rho e^{i\theta}) - i\rho e A_\mu e^{i\theta} = e^{i\theta}\partial_\mu\rho + \rho\partial_\mu(e^{i\theta}) - i\rho e A_\mu e^{i\theta}$$

$$= e^{i\theta}\partial_\mu\rho + e^{i\theta}i\rho\partial_\mu\theta - i\rho e A_\mu e^{i\theta} = e^{i\theta}[\partial_\mu\rho + i\rho(\partial_\mu\theta - eA_\mu)]$$

$$(D_\mu\varphi)^\dagger D_\mu\varphi = e^{-i\theta}e^{i\theta}[\partial_\mu\rho + i\rho(\partial_\mu\theta - eA_\mu)][\partial_\mu\rho - i\rho(\partial_\mu\theta - eA_\mu)]$$

$$= [\partial_\mu\rho + i\rho(\partial_\mu\theta - eA_\mu)][\partial_\mu\rho - i\rho(\partial_\mu\theta - eA_\mu)]$$

$$= (\partial_\mu\rho)^2 + \rho^2(\partial_\mu\theta - eA_\mu)^2$$

$$\varphi^\dagger\varphi = \mu^2\rho^2$$

$$\lambda(\varphi^\dagger\varphi)^2 = \lambda(\rho^2)^2 = \lambda\rho^4$$

$$\text{(H1)} \Rightarrow \mathcal{L} = -\tfrac{1}{4}F_{\mu\nu}F^{\mu\nu} + \rho^2(\partial_\mu\theta - eA_\mu)^2 + (\partial_\mu\rho)^2 + \mu^2\rho^2 - \lambda\rho^4 \qquad \text{(H2)}$$

Since $(\partial_\mu\rho)^2 = [\partial_\mu\tfrac{1}{\sqrt{2}}(v+\chi)]^2 = \tfrac{1}{2}[\partial_\mu(v+\chi)]^2 = \tfrac{1}{2}(\partial_\mu v + \partial_\mu\chi)^2$

$$= \tfrac{1}{2}(\partial_\mu\sqrt{\tfrac{\mu^2}{\lambda}} + \partial_\mu\chi)^2 = \tfrac{1}{2}(0 + \partial_\mu\chi)^2 = \tfrac{1}{2}(\partial_\mu\chi)^2$$

$$\rho^2(\partial_\mu - eA_\mu)^2 = \rho^2[e(A_\mu - \tfrac{1}{e}\partial_\mu)]^2 = \rho^2 e^2 B^2$$

$$\rho^2 = [\tfrac{1}{\sqrt{2}}(v+\chi)]^2 = \tfrac{1}{2}(v+\chi)^2 = \tfrac{1}{2}v^2 + v\chi + \tfrac{1}{2}\chi^2$$

$$\rho^2 e^2 B^2 = \left(\tfrac{1}{2}v^2 + v\chi + \tfrac{1}{2}\chi^2\right)e^2 B^2 = \tfrac{1}{2}v^2 e^2 B^2 + e^2 v\chi B^2 + \tfrac{1}{2}e^2\chi^2 B^2$$

$$= \tfrac{1}{2}M^2 B^2 + e^2 v\chi B^2 + \tfrac{1}{2}e^2\chi^2 B^2,$$

$$\text{(H2)} \Rightarrow \mathcal{L} = -\tfrac{1}{4}F_{\mu\nu}F^{\mu\nu} + \tfrac{1}{2}M^2 B^2 + e^2 v\chi B^2 + \tfrac{1}{2}e^2\chi^2 B^2 + \tfrac{1}{2}(\partial_\mu\chi)^2 + \mu^2\rho^2 - \lambda\rho^4 \qquad \text{(H3)}$$

Since $\mu^2\rho^2 = \mu^2(\tfrac{\chi+v}{\sqrt{2}})^2 = \tfrac{1}{2}\mu^2(\chi+v)^2$

$$= \tfrac{1}{2}\mu^2(\chi^2 + v^2 + 2\chi v) = \tfrac{1}{2}\mu^2\chi^2 + \tfrac{1}{2}\mu^2 v^2 + \mu^2\chi v$$

$$= \tfrac{1}{2}\mu^2\chi^2 + \tfrac{1}{2\lambda}\mu^4 + \mu^2\sqrt{\tfrac{\mu^2}{\lambda}}\chi = \tfrac{1}{2}\mu^2\chi^2 + \tfrac{1}{2\lambda}\mu^4 + \mu^3\chi\tfrac{1}{\sqrt{\lambda}}$$

$$\lambda\rho^4 = -\lambda(\frac{v+\chi}{\sqrt{2}})^4 = -\frac{\lambda}{4}(v+\chi)^4 = -\frac{\lambda}{4}[(v+\chi)^2]^2$$

$$= -\frac{\lambda}{4}(\chi^2 + v^2 + 2\chi v)^2 = -\frac{\lambda}{4}[(\chi^2 + v^2) + 2\chi v]^2$$

$$= -\frac{\lambda}{4}[(\chi^2 + v^2)^2 + 4\chi v(\chi^2 + v^2) + 4\chi^2 v^2]$$

$$= -\frac{\lambda}{4}(\chi^4 + v^4 + 2\chi^2 v^2 + 4v\chi^3 + 4\chi v^3 + 4\chi^2 v^2)$$

$$= -\frac{\mu^4}{4\lambda} - \frac{\lambda}{4}\chi^4 - \frac{1}{2}\mu^2\chi^2 - \mu^3\chi\frac{1}{\sqrt{\lambda}} - \sqrt{\lambda}\mu\chi^3 - \mu^2\chi^2$$

$$\mu^2\rho^2 - \lambda\rho^4 = \frac{1}{2}\mu^2\chi^2 + \frac{1}{2\lambda}\mu^4 + \mu^3\chi\frac{1}{\sqrt{\lambda}} - \frac{\mu^4}{4\lambda} - \frac{\lambda}{4}\chi^4 - \frac{1}{2}\mu^2\chi^2 - \mu^3\chi\frac{1}{\sqrt{\lambda}} - \sqrt{\lambda}\mu\chi^3 - \mu^2\chi^2$$

$$= \frac{\mu^4}{4\lambda} - \frac{\lambda}{4}\chi^4 - \sqrt{\lambda}\mu\chi^3 - \mu^2\chi^2$$

(H3) $\Rightarrow \quad \mathcal{L} = -\frac{1}{4}F_{\mu\nu}F^{\mu\nu} + \frac{1}{2}M^2 B^2 + e^2 v\chi B^2 + \frac{1}{2}e^2\chi^2 B^2$

$$+\frac{1}{2}(\partial_\mu\chi)^2 + \frac{\mu^4}{4\lambda} - \frac{\lambda}{4}\chi^4 - \sqrt{\lambda}\mu\chi^3 - \mu^2\chi^2 \tag{H4}$$

Notice the second term on the right side of the equation, the free gauge field $B^2$ is with a mass term $M^2$. H4 is a massive Lagrangian. Thus, for this derivation, we observe that $Der(L, L_m)$ holds; therefore, for $j = g[Der(L, L_m)], j \in X$. In other words, since its model $X$ is satisfiable, the mass predicate M is definable. QED.

*5. General discussion*

      Quantum field theory has successfully generated a number of general conclusions. It seems meaningful to disclose the logical forms of these conclusions. The formal representation of these conclusions would also enable us to address formal issues of these metalogical properties. Gödel numbering method seems useful as well as powerful to look at these issues from formal arithmetic perspectives. The results are obtained by integrating the twin theorems by Gödel and Tarski, which are briefly introduced in Appendix.

*Appendix     Twin theorems of Gödel and Tarski*

      We assume readers are familiar with the first order logic (**PL**), the first order theory ($\boldsymbol{\mathcal{N}}$), and the arithmetic theory (*N*). Intuitive natural numbers used in *N* are given by $n$, and the corresponding enumerers used in $\boldsymbol{\mathcal{N}}$ are denoted by bold $\boldsymbol{n}$. Enumerers are constructed by starting from the empty set ∅ and the so-called successor function, such that $\emptyset = \boldsymbol{0}$, { ∅} = **1**, {∅, {∅}} = **2**, and so forth. In the following, we introduce Godel's theorem first, and Tarski's theorem second.

*Definition of Expressibility:* If $R(a_1, \ldots, a_n)$ holds in $N$, then $P(\boldsymbol{a_1}, \ldots, \boldsymbol{a_n})$ is provable in $\mathcal{N}$. If $R(a_1, \ldots, a_n)$ does not hold in $N$, then $\neg P(\boldsymbol{a_1}, \ldots, \boldsymbol{a_n})$ is provable in $\mathcal{N}$.

*Definition of Consistency:* For any given formula "$L$" in $\mathcal{N}$, either $L$ is provable, or else $\neg L$ is provable, but not both.

For a given formula "$L$", denote its proof by $Bew(L)$. Assume $g(L) = i$ and $g(Bew(L)) = j$, where $i$ and $j$ are Gödel numbers. We introduce a relation $\mathcal{G}(i, j)$ in $\boldsymbol{N}$, and define a function term $G(\boldsymbol{i}, \boldsymbol{j})$ in $\mathcal{N}$. Gödel constructed a formula,

$$P(x) = \forall y \neg G(x, y),$$

in which $x$ is a free variable. Let $g[P(x)] = i$, by substituting $x$ with $i$, we can use,

$$S = P(\boldsymbol{i}) = \forall y \neg G(\boldsymbol{i}, y)$$

This is a so-called self-reflection sentence.

*Gödel First Theorem*   Neither $S$ nor $\neg S$ is provable in $\mathcal{N}$.

We now briefly sketch a proof. First, we prove that $S$ is not provable. Assume for contradiction that $S$ is provable, then it must have a proof, write $Bew(S)$, let $g(Bew(L)) = j$ and $g(S) = i$, so that $\mathcal{G}(i, j)$ in $\boldsymbol{N}$. By the expressibility, $G(\boldsymbol{i}, \boldsymbol{j})$ must be provable in $\mathcal{N}$; but $S$ said that for any $\boldsymbol{j}$, $\neg G(\boldsymbol{i}, \boldsymbol{j})$. This contradiction shows that the assumption is impossible. Hence, $S$ is not provable in $\mathcal{N}$.

Second, we prove that $\neg S$ is unprovable in $\mathcal{N}$. Assume for contradiction that $\neg S$ is provable. Then by consistency, $S$ is unprovable, so that for any $j$, $g(Bew(S)) \neq j$. Hence, for any $j$, $\mathcal{G}(i, j)$ does not hold in $\boldsymbol{N}$; by expressibility, $\neg G(\boldsymbol{i}, \boldsymbol{j})$, for any $j$. As such, by $\omega\_consistency$, we have $\forall y \neg G(\boldsymbol{i}, y)$, which means $S$ is provable in $\mathcal{N}$. This contradicts to the assumption that $S$ is unprovable. Thus $\neg S$ is unprovable in $\mathcal{N}$. QED- maybe define QED.

The above result shows that the consistency of $\mathcal{N}$ is independent of $\mathcal{N}$. Now let us speculate about what $S$ expresses. $S$ is a self-reflection sentence, it says that $S$ is unprovable, and we have just proved it above; thus, *$S$ is true, but not provable in $\mathcal{N}$*, which by definition means that $\mathcal{N}$ is incomplete. This is the well-known Gödel Incompleteness Theorem. We now turn to Tarski's indefinability theorem.

*Definition of Definability:* Let $g(P(x)) = m$, and $g(P(m)) = n$, we can hold a binary relation $d(m, n)$ in $\boldsymbol{N}$. Accordingly, we say $D(\boldsymbol{m}, \boldsymbol{n})$ is definable in $\mathcal{N}$, meaning $D(\boldsymbol{m}, \boldsymbol{n})$ has a model, which is not null.

Tarski introduced a new predicate of being true, denoted by $\boldsymbol{T}$, and he constructed a sentence below:

$$A(x) = \forall y[D(x, y) \rightarrow \neg T(y)]$$

Let $g[A(x)] = m$, substituting $x$ by $m$, we have

$$B(\boldsymbol{m}) = \forall y[D(\boldsymbol{m}, y) \to \neg T(y)]$$

Let $g[B(m)] = n$, we have $d(m, n)$ which holds in $\boldsymbol{N}$; hence, $D(x, y)$ is definable in $\boldsymbol{\mathcal{N}}$. Then, by standard logic, we can infer $\neg T(y)$. Now we show that $T$ is not definable, meaning its model is null. Denote $B(\boldsymbol{m})$ by $L$. If $L$ is pre-assumed as true, denote it by $Lt$, and write $Lg = g(L)$. As such, we may assume for contradiction that $T$ had a model $X$;

$$X = \{\, Lg \mid X \vDash L, \quad \text{i.e.,} Lt \text{ is presupposed to be true under } X\}$$

Since $Lg = g(Lt)$, i.e., $Lg$ is the Gödel number of $B(\boldsymbol{m})$. By the definition of $d(m, n)$, we have $d(m, Lg)$, hence $D(\boldsymbol{m}, \boldsymbol{Lg})$ is definable in $\boldsymbol{\mathcal{N}}$. However, recalling the logical structure of $B(\boldsymbol{m})$, which is a universally quantified conditional statement, we may infer $\neg T(\boldsymbol{Lg})$, i.e., $Lt$ is not true in model $X$; hence, $\boldsymbol{Lg} \notin X$, which shows that $X$ may only be null. In other words, since $T$ has no model, the truth predicate function is arithmetically undefinable.

*Acknowledgement* The author thanks Zhangju Liu for useful comments. I also thank my students who participated in my lectures and served as the first round proof readers.